\documentclass[12pt,a4paper]{article}
\usepackage[english]{babel}
\usepackage{amsmath}
\usepackage{amssymb}
\usepackage{graphicx}
\usepackage{color}

\newcommand{\nm}{\noalign{\smallskip}}
\newcommand{\ds}{\displaystyle}
\usepackage{pgf}
\usepackage{pgffor}
\usepgflibrary{plothandlers}
\usepackage{pgfplots}
\pgfplotsset{compat=newest}
 \pgfplotsset{width=15cm}
\pgfplotsset{plot coordinates/math parser=false}
\newlength\figureheight
\newlength\figurewidth

\newtheorem{de}{Definition}[section]
\newtheorem{theo}{Theorem}[section]

\newtheorem{prop}[theo]{Proposition}
\newtheorem{lemm}[theo]{Lemma}
\newtheorem{rem}[theo]{Remark}
\newcommand{\Proof}{\medskip\noindent {\bf Proof. }}

\newcommand{\cH}{{\cal H}}

\newcommand{\N}{\mathbb{N}}

\newcommand{\R}{\mathbb{R}}

\newcommand{\al}{\alpha}

\newcommand{\Om}{\Omega}

\newcommand{\ubf}{\mathbf{u}}
\newcommand{\hbf}{\mathbf{h}}

\newcommand{\nubf}{\boldsymbol{\nu}}
\newcommand{\epstilde}{\widetilde{\varepsilon}}

\newcommand{\piecewise}{\mathcal{C}_{S}^{k,\al}\big(\ol\Om\big)}

\newcommand{\bv}{\mathrm{BV}(\Om)}

\newcommand{\lone}{{L^{1}(\Om)}}

\newcommand{\cinf}{\mathcal{C}^\infty_0 \left( \Omega\right)}
\newcommand{\sbv}{\mathrm{SBV}(\Omega)}

\numberwithin{equation}{section} \numberwithin{figure}{section}

\newcommand{\ii}{\infty}

\newcommand{\nn}{\nonumber}

\newcommand{\dr}{\partial}

\newcommand{\g}{\nabla}

\newcommand{\ol}{\overline}
\newcommand{\su}{\subset}

\newcommand{\ra}{\longrightarrow}
\newcommand{\cqfd}{\hfill $\square$\\ \medskip}

\title{Mathematical modeling in full-field optical coherence elastography\thanks{\footnotesize This work was supported  by the
ERC Advanced
Grant Project MULTIMOD--267184.}}
\author{Habib Ammari\thanks{\footnotesize Department of Mathematics and Applications,
Ecole Normale Sup\'erieure, 45 Rue d'Ulm, 75005 Paris, France
(habib.ammari@ens.fr, pierre.millien@ens.fr,
laurent.seppecher@ens.fr).} \and Elie Bretin\thanks{\footnotesize
Institut Camille Jordan, Universit\'e de Lyon,  Lyon, F-69003, France
(bretin@cmap.polytechnique.fr).}
\and Pierre Millien\footnotemark[2] \and Laurent Seppecher\footnotemark[2] \and
Jin-Keun Seo\thanks{\footnotesize Department of Computational
Science and Engineering, Yonsei University, 50 Yonsei-Ro,
Seodaemun-Gu, Seoul 120-749, Korea (seoj@yonsei.ac.kr).}}

\begin{document}
\maketitle
\begin{abstract}
We provide a mathematical analysis of and a numerical framework for full-field optical coherence elastography, which has unique features including micron-scale resolution, real-time processing, and non-invasive imaging. We develop a novel algorithm for transforming volumetric optical images before and after the mechanical solicitation of a sample with sub-cellular resolution into quantitative shear modulus distributions. This has the potential to improve sensitivities and specificities in the biological and clinical applications of optical coherence tomography.
\end{abstract}

\bigskip

\noindent {\footnotesize Mathematics Subject Classification
(MSC2000): 35R30, 35B30.}


\section{Introduction}

Optical coherence tomography (OCT) is a non-invasive and a non-ionizing imaging technique that produces high-resolution images of biological tissues.  It performs optical slicing in the sample, to allow three-dimensional reconstructions of internal structures. Conventional optical coherence time-domain and frequency-domain tomographies require transverse scanning of the illumination spot in one or two directions to obtain cross-sectional or en face images, respectively.  Full-field OCT allows OCT to be performed without transverse scanning; the tomographic images are obtained by combining interferometric images acquired in parallel using an image sensor.  Both the transverse and the axial resolutions are of the order of $1\mu$m; see \cite{claude1,claude2}. We refer to \cite{otmar} for the mathematical modeling of OCT.

Elastography is an imaging-based technique for the estimation of the elastic properties of tissues. Given that the mechanical properties of tissues and cells are related to their structure and function, changes in those properties can reflect healthy or pathological states such as weakening of vessel walls or cirrhosis of the liver. Elastography can aid the identification of suspicious lesions, the diagnosis of various diseases and the monitoring of the effectiveness of  treatments (see \cite{elasto1,elasto2}). Different imaging modalities (e.g., ultrasound and magnetic resonance imaging) can be used to measure tissue displacements and  to estimate the resulting tissue stiffness and viscosity. Magnetic resonance elastography  is relatively expensive, due to the high magnetic field environment, which requires specifically designed equipment.
Several reconstruction approaches for elastography have been derived
\cite{ammari2008introduction, ammarielasto, garapon, seobook}.

In \cite{nahas20133d}, elastographic contrast has been combined with full-field OCT with the aim of creating a virtual palpation map at the micrometer scale.
The idea is to register a volumetric optical image before and after mechanical solicitation
of the sample. Based on the assumption that the density of the optical scatterers is advected by the deformation, the displacement map can be first estimated. Then, using a quasi-incompressible model for the tissue elasticity, the shear modulus distribution can be reconstructed from the estimated displacement map.

The OCT elastography is able to perform displacement measurements with sub-cellular resolution.  It enables a more precise characterization of tissues than that achieved using ultrasound or magnetic resonance elastography;  therefore, it provides a more accurate assessment of microscale variations of elastic properties. A map of mechanical properties added as a supplementary contrast mechanism to morphological images could aid diagnosis. The technique costs less than other elastography techniques.

The mapping of mechanical properties was first introduced  to OCT imaging by Schmitt \cite{schmitt}, who measured displacements as small as a few micrometers in heterogeneous gelatin  phantoms containing scattering particles in addition to living skin. Various subsequent applications have employed  OCT methods in elastography; these include  dynamic and full-field optical coherence elastography (see \cite{dynamic, elasto, cardiac}).

In all of the aforementioned techniques, transforming the OCT images before and after the application of a load into quantitative maps of the shear modulus is a challenging problem.

In this paper we present a mathematical and numerical framework for the OCT-elastography experiment described in \cite{nahas20133d}. Using the set of images before and after mechanical solicitation we design a novel method to reconstruct the shear modulus distribution inside the sample.

To mathematically formulate the problem, let $\Omega_0 \su \R^d, d=2,3,$ and let $\varepsilon_0$ be the known  piecewise smooth optical index of the medium, and $\mu$ be its  shear modulus.
In this paper we consider heterogeneous (unknown) shear modulus distributions.
The medium is solicited mechanically. Since compression modulus of biological media is four order of magnitude larger than the shear modulus, it can be shown that
the displacement map $\ubf$ obeys the linearized equations of incompressible fluids or the
Stokes system \cite{ammari2008introduction, ammarielasto,garapon}. The model problem is then the following Stokes system in a heteregeneous medium which reads:
\begin{equation}\label{eqmu}\left\{\begin{aligned}
\g\cdot \left(\mu ( \g \ubf + \g \ubf^T) \right)+ \g p=0 \quad &  \mbox{in } \ \Omega_0,\\
\g\cdot \ubf = 0 \quad & \mbox{in } \ \Omega_0,\\
\ubf =\mathbf{f} \quad & \mbox{on } \ \dr \Omega_0,
\end{aligned}\right.
\end{equation}
where superposed $T$ denotes the transpose and the real-valued vector $\mathbf{f} $ satisfies the compatibility condition $\int_{\partial \Omega_0} \mathbf{f}
\cdot {\nubf} =0$ with ${\nubf}$ being the outward normal at $\partial \Omega_0$.

Throughout this paper, we assume that  $\mu \in \mathcal{C}^{0,1}(\overline \Omega_0)$
and $\mathbf{f} \in \mathcal{C}^{2}(\partial \Omega_0)^d$.  From \cite{chen, giaquinta, yanyan}, (\ref{eqmu}) has a unique solution
$\ubf \in \mathcal{C}^{1}(\overline{\Omega_0})^d$ .
Moreover,  there exists a positive constant $C$ depending only on $\mu$ and $\Omega_0$ such that
$$
|| \ubf ||_{\mathcal{C}^{1}(\overline{\Omega_0})^d} \leq C || {\mathbf{f}}||_{\mathcal{C}^{2}({\partial \Omega_0})^d}.$$

Using a second OCT scan, one has access to the optical index of the deformed medium $\varepsilon_u(\widetilde{\mathbf{x}}),\ \forall\; \widetilde{\mathbf{x}}\in \Omega_u$, where
$\Omega_u$ is defined by
$$\Omega_u=\{\mathbf{x} +\ubf(\mathbf{x}), \ \mathbf{x} \in \Omega_0\}.$$

The new optical index is linked to the original one by
\begin{equation}
\varepsilon(\mathbf{x})=\varepsilon_u\left(\mathbf{x}+\ubf(\mathbf{x})\right),\quad \forall\; \mathbf{x} \in \Omega_0.
\end{equation}
The goal is to reconstruct the shear modulus map $\mu$ on $\Omega_0$ from the functions $\varepsilon$ and $\varepsilon_u$. We first prove that, in two dimensions, if the direction of $\displaystyle{\frac{\g \varepsilon}{\vert \g \varepsilon \vert }}$ is not constant in a neighborhood of $\mathbf{x}$, then the displacement field $\ubf$ at $\mathbf{x}$ can be approximately reconstructed. In three dimensions, one shall assume that the vectors $\displaystyle{\frac{\g \varepsilon (\mathbf{y})}{\vert \g \varepsilon (\mathbf{y})\vert }}$ are not coplanar for $\mathbf{y}$ a neighborhood of $\mathbf{x}$.
Hence, the reconstructed value of $\ubf(\mathbf{x})$ serves as an initial guess for  the minimization of the discrepancy between computed and measured changes in the optical index. Then, we compute an element of the subgradient \cite{clarke}
of the discrepancy functional. Finally, we implement a minimization scheme to retrieve the shear modulus map from the reconstructed displacements.

The paper is organized as follows.
Section \ref{sect2} is devoted to some mathematical preliminaries. In section \ref{sect3} we consider piecewise smooth $\varepsilon$ functions and first derive a leading-order Taylor expansion of $\varepsilon_u$ as $|| \ubf||_{\mathcal{C}^1}$ goes to zero. Then we provide an initial guess by linearization. Finally, we prove the Fr\'echet differentiability of the discrepancy functional between the measured and the computed advected images. The displacement field inside the sample can be obtained as the minimizer of such functional.
Section \ref{sect4} is devoted to the reconstruction of the shear modulus from the displacement measurements. In section \ref{sect5} we present some numerical results to highlight the viability and the performance of the proposed algorithm. The paper ends with a short discussion.

\section{Preliminaries} \label{sect2}

Let $\Omega$ be a bounded smooth domain in $\mathbb{R}^d, d=2,3$. We start by defining a class of piecewise smooth functions.
\begin{de}\label{defpiecewise} For any $k\in\N$, $\al\in ]0,1[$, for any curve $S$ of class $\mathcal{C}^{1,\alpha}$ for some $0<\alpha<1$ such that  $\Om\setminus S$ is a union of connected domains $\Om_i, i=1,2,\cdots n$,   we define $\piecewise$ to be the class of functions $f:\Om\ra\R$ satisfying

\begin{equation}\begin{aligned}
 f|_{\Om_i}\in \mathcal{C}_{S}^{k,\al}\big(\ol\Om_i\big) \quad \forall \; i=1,\cdots n.
\end{aligned}\end{equation}
\end{de}

\begin{de} We define $\bv$ as the subspace of $\lone$ of all the
functions $f$ whose  weak derivative $Df$ is a finite Radon measure.
In other terms, $f$ satisfies

\begin{equation}\nn\begin{aligned}
\int_\Omega f\g\cdot \mathbf{F} \leq C \sup_{x\in\Omega}|\mathbf{F}|,\quad \forall \; \mathbf{F}\in \mathcal{C}^1_0(\Omega)^d
\end{aligned}\end{equation}
\end{de}
for some positive constant $C$ with $ \mathcal{C}^1_0(\Omega)$ being the set of compactly supported $\mathcal{C}^1$ functions.

The derivative of a function $f\in\bv$ can be decomposed as

\begin{equation}\nn\begin{aligned}
Df=\g f\mathcal{H}^d+[f]\nubf_s\cH^{d-1}_S+D_cf,
\end{aligned}\end{equation}
where $\mathcal{H}^d$ is the Lebesgue measure on $\Om$, $\mathcal{H}_S^{d-1}$ is the
surface Hausdorff measure on a rectifiable surface $S$, $\nubf_S$ is a
normal vector defined a.e. on $S$, $\g f\in\lone$ is the smooth
derivative of $f$, $[f]\in L^1(S,\mathcal{H}^{d-1}_S)$ is the jump of $f$ across $S$ and
$D_c f$ is a vector measure supported on a set of Hausdorff dimension less than $(d-1)$, which means that its ${(d-1)}$-Hausdorff-measure is zero; see \cite{alberti_amb}.

\begin{de} We define $\sbv$ as the subspace of $\bv$ of all the functions
$f$ satisfying $D_c f=0$. \end{de}

\begin{de} For any $1\leq p\leq +\ii$, we define
\begin{equation}\nn\begin{aligned}
\mathrm{SBV}^p(\Om)=\left\{f\in\sbv\cap L^p(\Omega),\ \g f\in L^p(\Omega)^d\right\}.
\end{aligned}\end{equation}
\end{de}

As $\sbv$ is a good model for piecewise-$W^{1,1}$ functions, the space
$\mathrm{SBV}^p(\Om)$ can be seen as the space of piecewise-$W^{1,p}$ functions. Here,
$W^{1,p}(\Omega)=\{ f \in L^p(\Omega),\ \g f\in L^p(\Omega)^d \}$ for $p\geq 1$.

Note that the space $\mathrm{SBV}^\ii(\Om)$ is a nice definition of piecewise Lipschitz
function. Note also that $\piecewise \subset \mathrm{SBV}^p(\Omega)$.

From now on, we assume that the optical index in the medium $\varepsilon$ belongs to $\piecewise$, which is a simple but good model for a discontinuous medium. Some of the following propositions are true for more general maps $\varepsilon\in \mathrm{SBV}(\Omega)$. In these propositions we only assume that $\varepsilon$ is in $\mathrm{SBV}(\Omega)$.
\section{Displacement field measurements} \label{sect3}
\subsection{First order approximation}

Let  $\Omega\Subset \left(\Omega_0\cap \Omega_u\right)$ be a smooth simply connected domain.
On $\Omega$, we have
\begin{equation*}
\begin{aligned}\varepsilon_u&=&\varepsilon \circ \left( \mathbb{I} + \ubf \right)^{-1}\\ \varepsilon&=&\varepsilon_u \circ \left( \mathbb{I} + \ubf \right) ,\end{aligned}
\end{equation*}
where $\mathbb{I}$ is the $d\times d$ identity matrix. 

\begin{prop}\label{firstorder}Let $\varepsilon \in \mathrm{BV}(\Omega)$ and let $\ubf \in \mathcal{C}^{1}(\overline{\Omega})^d$ be such that $\Vert \ubf \Vert_{\mathcal{C}^1(\overline{\Omega})^d} <1$.  Then, for any $\psi \in \mathcal{C}^1_0(\Omega)$, we have
\begin{equation} \label{estcv}
\left\vert\int_\Omega \left( \varepsilon - \varepsilon_u\right) \psi - \int_\Omega \psi\ubf \cdot  D\varepsilon \right\vert\leq C \Vert \ubf\Vert_{\mathcal{C}^0(\overline{\Omega})^d}\Vert \ubf \Vert_{\mathcal{C}^1(\overline{\Omega})^d} \Vert \psi \Vert_{\mathcal{C}^1_0(\Omega)} \vert  \varepsilon \vert_{\mathrm{TV}(\Omega)},
\end{equation}
where the constant $C$ is independent of $\psi$ and $\vert\; \vert_{\mathrm{TV}(\Om)}$ denotes the total variation semi-norm.
Estimate (\ref{estcv}) yields that $\displaystyle{\frac{\varepsilon_u - \varepsilon + \ubf \cdot D\varepsilon}{\Vert \ubf \Vert_{\mathcal{C}^0(\overline{\Omega})^d}}}$ weakly converges  to $0$ in $\mathcal{C}^1_0(\Omega)$ when $\Vert \ubf \Vert_{\mathcal{C}^1(\overline{\Omega})^d}$ goes to $0$.
\end{prop}
\Proof For each $t\in [0,1]$, define $ \phi_t$ by $\phi_t^{-1}(\mathbf{x})= \mathbf{x} +t \ubf(\mathbf{x})$.
Let $\eta >0$ be a small parameter, and $\varepsilon^{(\eta)}$ be a smooth function such that $\Vert\varepsilon-\varepsilon^{(\eta)}\Vert_{L^1(\Omega)} \rightarrow 0$, and $\vert \varepsilon^{(\eta)} \vert_{\mathrm{TV}(\Omega)} \rightarrow \vert \varepsilon \vert_{\mathrm{TV}(\Omega)}$ as $\eta \rightarrow 0$. Analogously, we define $\varepsilon^{(\eta)}_u$ to be the smooth approximation of
$\varepsilon_u$ given by
$$
\varepsilon^{(\eta)}_u (\mathbf{x})= \varepsilon^{(\eta)} \circ \phi_1(\mathbf{x}).
$$
From
\begin{equation*}
\varepsilon^{(\eta)}_u(\mathbf{x}) - \varepsilon^{(\eta)}(\mathbf{x}) = \left(\varepsilon^{(\eta)}\circ \phi_1\right)(\mathbf{x}) - \left(\varepsilon^{(\eta)}\circ \phi_0\right)(\mathbf{x}),
\quad \forall \; \mathbf{x} \in \Omega,
\end{equation*}
we have
\begin{equation*}
\varepsilon^{(\eta)}_u(\mathbf{x}) - \varepsilon^{(\eta)}(\mathbf{x}) = \int_0^1 \g \varepsilon^{(\eta)}(\phi_t(\mathbf{x}))\cdot \dr_t \phi_t(\mathbf{x}) dt, \quad \forall \; \mathbf{x}\in \Omega.
\end{equation*}
Therefore, for $\psi \in \cinf$ with $\cinf$  being the set of compactly supported $\mathcal{C}^\infty$ functions, 
\begin{multline}
\int_\Omega \left[\varepsilon^{(\eta)}_u(\mathbf{x}) - \varepsilon^{(\eta)}(\mathbf{x}) + \g \varepsilon^{(\eta)}(\mathbf{x}) \cdot \ubf(x)\right] \psi(\mathbf{x}) d\mathbf{x} =\\  \int_\Omega \left[\int_0^1 \g \varepsilon^{(\eta)}(\phi_t(\mathbf{x}))\cdot \dr_t \phi_t(\mathbf{x}) dt \right]\psi(\mathbf{x})  d\mathbf{x} + \int_\Omega \g\varepsilon^{(\eta)}(\mathbf{x}) \cdot \ubf(\mathbf{x}) \psi(\mathbf{x})d\mathbf{x},  \quad \forall \; \mathbf{x}\in \Omega.
\end{multline}
By a change of variables in the first integral and using the fact that $$\dr_t \phi_t(\mathbf{x}) = -\dr_{\mathbf{x}} \phi_t(\mathbf{x}) \dr_t \phi_t^{-1}(\mathbf{y})|_{\mathbf{y}=\phi_t(\mathbf{x})},$$ we get, for all $\mathbf{x}\in \Omega$,
\begin{multline*}
\int_0^1 \left[\int_\Omega  \g \varepsilon^{(\eta)}(\phi_t(\mathbf{x}))\cdot \dr_t \phi_t(\mathbf{x}) \psi(\mathbf{x})d\mathbf{x} \right]dt = \\   -\int_0^1 \int_{\Omega } \g \varepsilon^{(\eta)} (\mathbf{y})  \cdot \left[\dr_{\mathbf{x}} \phi_t(\phi_t^{-1}(\mathbf{y})) \dr_t \phi_t^{-1}(\mathbf{y}) \right] \vert \text{det } \dr_{\mathbf{x}} \phi_t^{-1}(\mathbf{y})\vert \psi \left(\phi_t^{-1}(\mathbf{y})\right)d\mathbf{y} dt .
\end{multline*}
Here, $\text{det }$ denotes the determinant of a matrix. 
Since $$\forall \; (\mathbf{y},t)\in \Omega \times [0,1],  \quad
\dr_t \phi_t^{-1}(\mathbf{y}) = \ubf(\mathbf{y}), $$  $$\dr_{\mathbf{y}} \phi_t^{-1}(\mathbf{y}) = \mathbb{I} + t \nabla \ubf (\mathbf{y}),$$  and $$\dr_{\mathbf{x}} \phi_t(\phi_t^{-1}(\mathbf{y})) \dr_{\mathbf{y}} \phi_t^{-1}(\mathbf{y}) = \mathbb{I},$$
we can write
\begin{multline*}
\int_0^1 \int_\Omega \left[ \g \varepsilon^{(\eta)}(\phi_t(\mathbf{x})) \cdot \dr_t \phi_t(\mathbf{x}) \psi(\mathbf{x})d\mathbf{x} \right]dt  =\\ -\int_0^1 \int_\Omega \g\varepsilon^{(\eta)}(\mathbf{y}) \cdot \left[ \left(\mathbb{I} + t \nabla \ubf (\mathbf{y})\right)^{-1} \ubf(\mathbf{y}) \right] \vert \text{det } \mathbb{I} + t \nabla \ubf (\mathbf{y}) \vert \psi\left(\phi_t^{-1}(\mathbf{y})\right) d\mathbf{y} dt,
\end{multline*}
and hence,
\begin{multline}\label{expression}
\int_\Omega \left[\varepsilon_u^{(\eta)}(\mathbf{x}) - \varepsilon^{(\eta)}(\mathbf{x}) + \g \varepsilon^{(\eta)}(\mathbf{x}) \cdot \ubf(\mathbf{x})\right] \psi(\mathbf{x}) d\mathbf{x} = \\ \int_0^1 \int_\Omega \g \varepsilon^{(\eta)} (\mathbf{x}) \cdot   \ubf(\mathbf{x}) \big[  \psi(\mathbf{x})-\psi\left(\phi_t^{-1}(\mathbf{x}) \right)  \big] d\mathbf{x} dt
\\ + \int_0^1 \int_\Omega \g \varepsilon^{(\eta)} (\mathbf{x}) \cdot \left( \left[ \left(\mathbb{I} + t \nabla \ubf (\mathbf{x})\right)^{-1} \vert \text{det } \mathbb{I} + t\nabla \ubf (\mathbf{x}) \vert - \mathbb{I} \right] \ubf(\mathbf{x})\right) \psi\left(\phi_t^{-1}(\mathbf{x}) \right) d \mathbf{x} dt.
\end{multline}
The first term in the right-hand side of (\ref{expression}) can be estimated as follows:
\begin{equation*}
\left\vert \int_0^1 \int_\Omega \g \varepsilon^{(\eta)} (\mathbf{x}) \cdot   \ubf(\mathbf{x}) \big[  \psi(\mathbf{x})-\psi\left(\phi_t^{-1}(\mathbf{x}) \right)  \big] d\mathbf{x} dt \right\vert \leq  \Vert \ubf \Vert_{\mathcal{C}^0(\overline{\Omega})^d}^2 \Vert \g \varepsilon^{(\eta)} \Vert_{L^1(\Omega)^d} \Vert \g \psi \Vert_{\mathcal{C}^0({\Omega})^d}.
\end{equation*}
Let $\text{tr }$ denote the trace of a matrix. Using the fact that
\begin{equation*}
\left(\mathbb{I} + t \nabla \ubf \right)^{-1}=\sum_{i=0} (-1)^i \left( t \nabla \ubf \right)^i,
\end{equation*}
which follows from $|| \ubf||_{\mathcal{C}^1(\Omega)^d} < 1$, and
\begin{equation*}
\text{det } \left(\mathbb{I} + t\nabla \ubf \right) =\left\{
\begin{aligned}
&1-\text{tr } t \nabla \ubf + \text{det } t \nabla \ubf & \text{ if } d=2, \\
& 1+ \text{tr } t \nabla \ubf - \frac{1}{2}\left[ \left( \text{tr }t \nabla \ubf\right)^2 - \text{tr } \left(t \nabla \ubf\right)^2 \right]  + \text{det } t \nabla \ubf & \text{ if } d=3,
\end{aligned}\right.
\end{equation*}
we get
\begin{multline*}
\int_0^1 \int_\Omega \g \varepsilon^{(\eta)} (\mathbf{x}) \cdot \ubf(\mathbf{x}) \left[ \left(\mathbb{I} + t \g \ubf (\mathbf{x})\right)^{-1} \vert \text{det } \mathbb{I} + t\nabla \ubf (\mathbf{x}) \vert - \mathbb{I} \right]  \psi\left(\phi_t^{-1}(\mathbf{x}) \right) dxdt \\
\leq \Vert \ubf \Vert_{\mathcal{C}^0(\overline{\Omega})^d}
\Vert \ubf \Vert_{\mathcal{C}^1(\overline{\Omega})^d} \Vert \g \varepsilon^{(\eta)} \Vert_{L^1(\Omega)^d}  \Vert\psi\Vert_{\mathcal{C}^0({\Omega})},
\end{multline*}
which is the desired estimate for the second term in the right-hand side of (\ref{expression}).

Now, we can deduce the final result by density when $\eta \rightarrow 0$.
Since $\ubf \in \mathcal{C}^1(\Omega)^d$ and $\psi \in \mathcal{C}_0^1(\Omega)$, we can write $$\int_\Omega \psi \ubf \cdot \g\varepsilon^{(\eta)} =-\int_\Omega \g\cdot (\psi \ubf) \varepsilon^{(\eta)}. $$
Since $\Vert \varepsilon^{(\eta)} -\varepsilon\Vert_{L^1(\Omega)} \rightarrow 0$, we have $$\int_\Omega \g \cdot(\psi \ubf) \varepsilon^{(\eta)} \rightarrow \int_\Omega \g\cdot (\psi \ubf) \varepsilon.$$ As $\vert \varepsilon^{(\eta)}\vert_{\mathrm{TV}(\Omega)} \rightarrow \vert \varepsilon \vert_{\mathrm{TV}(\Omega)}$, we arrive at (\ref{estcv}) and the proof of the proposition is complete. \cqfd

\subsection{Local recovery via linearization}
Assuming that $\varepsilon \in \mathrm{SBV}^2(\Omega)$, we can write
 \begin{equation*}
D\varepsilon= \g \varepsilon \mathcal{H}^d + [\varepsilon]_S\nubf_S \mathcal{H}^{d-1}_S,
\end{equation*}
where $\nubf_S$ is the outward normal at the oriented surface $S$ of discontinuity  of $\varepsilon$.

The data consists of $\varepsilon$ and $\varepsilon_u$ on $\Omega$. In order to reconstruct $\ubf$, we can use the first order approximation of $\varepsilon-\varepsilon_u$: $$\varepsilon-\varepsilon_u\approx \ubf\cdot D\varepsilon,$$  given by Proposition \ref{firstorder}. These data can be decomposed into two parts:
\begin{equation*}
\ubf\cdot D\varepsilon(\cdot) = \ubf \cdot \g \varepsilon \mathcal{H}^d + [\varepsilon]_S \ubf \cdot \nubf_S \mathcal{H}^{d-1}_S = d_{\mathrm{reg}}\mathcal{H}^d + d_{\mathrm{sing}} \mathcal{H}^{d-1}_S.
\end{equation*}

Let $w$ be a mollifier supported on $[-1,1]$. For any $\delta >0$, we define $$w_\delta=\frac{1}{\delta^d} w\left(\frac{\cdot}{\delta}\right), $$ and  introduce $$\ubf_\delta(\mathbf{x})=\int_{\Omega} \ubf(\mathbf{y}) w_\delta(\vert \mathbf{x} -\mathbf{y} \vert) d\mathbf{y}.$$ Since $\ubf$ is smooth, for any $\mathbf{x} \in \Omega$, $\ubf_\delta(\mathbf{x})$ is a good approximation of $\ubf$ on the ball with center $\mathbf{x}$ and radius $\delta$.

We want to find an approximate value for $\ubf_\delta$ from the optical measurements and use it as an initial guess in an optimization procedure. For doing so, we introduce the functional $J_{\mathbf{x}} : \R^d \longrightarrow \mathbb{R}$ given by
\begin{multline*}
  \ubf \longmapsto J_{\mathbf{x}}(\ubf)= \int_{\Omega} \vert \g \varepsilon(\mathbf{y}) \cdot \ubf - d_{\mathrm{reg}}(\mathbf{y})\vert^2 w_\delta(\vert \mathbf{x} -\mathbf{y} \vert) d\mathbf{y} \\+  \int_{\Omega} \vert [\varepsilon]_S \ubf \cdot \nubf_S  - d_{\mathrm{sing}}(\mathbf{y}) \vert^2w_\delta(\vert \mathbf{x}-\mathbf{y}\vert) d\mathbf{y},
\end{multline*}
and look for minimizers of $J_{\mathbf{x}}$ in $\R^d$. The gradient of $J_{\mathbf{x}}$ can be explicitly computed as follows:
\begin{multline*}
\g J_{\mathbf{x}}(\ubf) = 2\int_{\Omega} \left(\g\varepsilon(y) \cdot \ubf - d_{\mathrm{reg}}(\mathbf{y})\right) \g \varepsilon(\mathbf{y}) w_\delta(\vert \mathbf{x}-\mathbf{y}\vert) d\mathbf{y} \\+ 2\int_{\Omega} \left( [\varepsilon]_S(\mathbf{y}) \ubf \cdot \nubf(\mathbf{y}) - d_{\mathrm{sing}}(\mathbf{y}) \right)[\varepsilon]_S(\mathbf{y})\nubf(\mathbf{y})w_\delta(\vert \mathbf{x} -\mathbf{y} \vert) d\mathbf{y}.
\end{multline*}

In the case where $\varepsilon$ has no jumps, $J_{\mathbf{x}}$ is a quadratic functional and we have
\begin{equation}
\g J_{\mathbf{x}}(\ubf) =0 \Leftrightarrow \ubf^T \left(\int_{\Omega} w_\delta(\vert \mathbf{x}-\mathbf{y}\vert)  \g \varepsilon(\mathbf{y}) \g \varepsilon^T (\mathbf{y})  d\mathbf{y}  \right)= \int_{\mathbf{x} +\delta B} d_{\mathrm{reg}}(\mathbf{y}) w_\delta(\vert \mathbf{x}-\mathbf{y} \vert) \g \varepsilon(\mathbf{y}) d\mathbf{y},
\end{equation}
where $B$ is the ball with center $0$ and radius $1$.

If the matrix $\displaystyle{\int_{\Omega}w_\delta(\vert \mathbf{x}-\mathbf{y}\vert)  \g \varepsilon(\mathbf{y}) \g \varepsilon^T (\mathbf{y})}$ is invertible, then the minimizer is given by
\begin{equation} \label{estinitial}
\ubf^T  =  \left(\int_{\Omega}  w_\delta(\vert \mathbf{x}-\mathbf{y} \vert)  \g \varepsilon(\mathbf{y}) \g \varepsilon^T (\mathbf{y}) d\mathbf{y} \right)^{-1}\int_{\mathbf{x} +\delta B} d_{\mathrm{reg}}w_\delta(\vert \mathbf{x}-\mathbf{y}\vert)  \g \varepsilon(\mathbf{y})d\mathbf{y}.
\end{equation}
The following proposition gives a sufficient condition for the invertibilty of the matrix $\displaystyle{\int_{\Omega}w_\delta(\vert \mathbf{x}-\mathbf{y}\vert)  \g \varepsilon(\mathbf{y}) \g \varepsilon^T (\mathbf{y})}$.
\begin{prop}\label{prop1} Suppose that $\varepsilon$ has no jumps and $d=2$. Assume $\mathbf{x}+\delta B\subset \Om$.
Then, if all vectors $\g \varepsilon$ in $\{ \mathbf{y}~:~ w_\delta(|\mathbf{y}- \mathbf{x}|)\neq 0\}$ are not  collinear, then the matrix $$\displaystyle{\int_{\Omega} w_\delta(\vert \mathbf{x}-\mathbf{y} \vert)  \g \varepsilon(\mathbf{y}) \g \varepsilon^T (\mathbf{y})d\mathbf{y} }$$ is invertible.
\end{prop}

\Proof  Writing
$$\forall \; \mathbf{y} \in \mathbf{x} +\delta B, \quad \g \varepsilon (\mathbf{y}) = u(\mathbf{y}) \mathbf{e}_1 + v(\mathbf{y}) \mathbf{e}_2,$$ where $\{\mathbf{e}_1,\mathbf{e}_2\}$ is the cannonical basis of $\mathbb{R}^2$,  it follows that \begin{equation*}
\g \varepsilon \g \varepsilon^T (\mathbf{y})= u^2(\mathbf{y})  \mathbf{e}_1 \mathbf{e}_1^T + v^2(\mathbf{y}) \mathbf{e}_2 \mathbf{e}_2^T +u(\mathbf{y})v(\mathbf{y})\left( \mathbf{e}_1
\mathbf{e}_2^T + \mathbf{e}_2 \mathbf{e}_1^T \right), \quad \forall \; \mathbf{y} \in \mathbf{x}+\delta B.
\end{equation*}
Computing the convolution with respect to $w_\delta$, we get
\begin{multline*}
\int_{\Omega} w_\delta(\vert \mathbf{y} - \mathbf{y}\vert) \g \varepsilon (\mathbf{y})\g \varepsilon^T (\mathbf{y}) d\mathbf{y} = \left(\int_{\Omega}  u^2(\mathbf{y}) w_\delta(\vert \mathbf{y} - \mathbf{x}\vert)d\mathbf{y} \right) \mathbf{e}_1 \mathbf{e}_1^T\\ + \left(\int_{\Omega}  v^2(\mathbf{y}) w_\delta(\vert \mathbf{y} - \mathbf{x}\vert)d\mathbf{y}\right) \mathbf{e}_2 \mathbf{e}_2^T + \left(\int_{\Omega}  u(\mathbf{y}) v(\mathbf{y})w_\delta^T(\vert \mathbf{y}- \mathbf{x}\vert) d\mathbf{y} \right) \left(\mathbf{e}_1 \mathbf{e}_2^T + \mathbf{e}_2 \mathbf{e}_1^T \right).
\end{multline*}
This matrix is not invertible if and only if
\begin{equation*}
\left(\int_{\Omega}  u^2(\mathbf{y}) w_\delta(\vert \mathbf{y} - \mathbf{x}\vert)d\mathbf{y}\right) \left(\int_{\Omega}  v^2(\mathbf{y}) w_\delta(\vert \mathbf{y} - \mathbf{x}\vert)d\mathbf{y} \right) = \left(\int_{\Omega}  u(\mathbf{y})v(\mathbf{y}) w_\delta(\vert \mathbf{y} - \mathbf{x}\vert)d\mathbf{y}\right)^2,
\end{equation*} which is exactly the equality case in weighted Cauchy-Schwarz inequality. So,  if there exist two points $ \mathbf{y_1}, \mathbf{y_2}\in \{ \mathbf{y}:~ w_\delta(|\mathbf{y}- \mathbf{x}|)\neq 0\}$ such that $\g \varepsilon(\mathbf{y_1})\times \g \varepsilon(\mathbf{y_2})\neq 0$, then $u$ is not proportional to $v$, and the matrix is invertible.\cqfd
\begin{rem} Assuming that $\g \varepsilon (\mathbf{y})\neq 0$ for  $\mathbf{y}\in \mathbf{x}+\delta B\subset \Om$, Proposition \ref{prop1} gives that the direction of
 $\displaystyle{\frac{\g \varepsilon}{\vert \g \varepsilon \vert }}$ in not constant in  $\mathbf{x}+\delta B\subset \Om$ if and only if
 $$\displaystyle{\int_{\mathbf{x}+\delta B}  \g \varepsilon(\mathbf{y}) \g \varepsilon^T (\mathbf{y})d\mathbf{y} } \quad\mbox{is invertible.}$$
Hence, under the above condition on $\varepsilon $ in the neighborhood $\mathbf{x}+\delta B$, the displacement field $\ubf$ at $\mathbf{x}$ can be approximately reconstructed.
\end{rem}
\begin{rem}
By exactly the same arguments as those in two dimensions, one can prove that in the three-dimensional case, if  all vectors $\g \varepsilon$ in $\{ \mathbf{y}~:~ w_\delta(|\mathbf{y}- \mathbf{x}|)\neq 0\}$ are not coplanar, then the matrix $$\displaystyle{\int_{\Omega} w_\delta(\vert \mathbf{x}-\mathbf{y} \vert) \g \varepsilon(\mathbf{y}) \g \varepsilon^T (\mathbf{y})d\mathbf{y} }$$ is invertible.

On the other hand, in the  case where $\varepsilon$ is  piecewise smooth, one can first detect the surface of jumps of $\varepsilon$ using for example an edge detection algorithm \cite{canny, canny2} and then apply the proposed local algorithm in order to have a good approximation of $\ubf$ in the domains where $\varepsilon$ is smooth.
\end{rem}

\subsection{Minimization of the discrepancy functional}

Let $\varepsilon \in \piecewise$, where $S$ is the surface of discontinuity. For the sake of simplicity we assume that $\Omega \setminus S$ is the union of two connected domains $\Omega_1 \cup \Omega_2$. Therefore, $\varepsilon$ can be written as \begin{equation}\label{hypothese epsilon}\varepsilon=\varepsilon_1 \chi_{\Omega_1} + \varepsilon_2 \chi_{\Omega_2}
\end{equation} with $\varepsilon_i \in \mathcal{C}^1(\ol \Omega_i)$, for $i=1,2$.

Denote $\ubf^*$ the applied (true) displacement on $\Omega$ (as defined in (\ref{eqmu})) and  $\epstilde$ the measured deformed optical index given by \begin{equation*}
\epstilde = \varepsilon \circ \left( \mathbb{I} + \ubf^*\right)^{-1}.
\end{equation*}
The following result holds.
\begin{prop}\label{diffentiabilite}
Let $\varepsilon$ verify (\ref{hypothese epsilon}), $\ubf^* \in \mathcal{C}^{1}(\Omega)^d$ be the solution of (\ref{eqmu}), and $\epstilde = \varepsilon \circ \left(\mathbb{I}+ \ubf^*\right)^{-1}$.  Suppose that $\Omega_2\Subset \Omega$. Then, the functional $I$ defined by
\begin{equation}
\begin{aligned} I:
\mathcal{C}^1(\Omega)^d &\longrightarrow \mathbb{R}, \\
  \ubf &\longmapsto I(\ubf)= \int_{\Omega} \vert \widetilde{\varepsilon} \circ (\mathbb{I}+\ubf) - \varepsilon \vert^2\, d\mathbf{x}
\end{aligned}
\end{equation} has a nonempty subgradient. Let $\boldsymbol{\xi}$ in the dual of $\mathcal{C}^1(\Omega)^d$ be given by
\begin{equation} \label{defxi}
\boldsymbol{\xi}: \mathbf{h} \mapsto 2\int_{\Omega} [\epstilde(\mathbf{x}+\ubf) -\varepsilon(\mathbf{x})] \hbf(\mathbf{x})\cdot D\epstilde\circ(\mathbb{I}+\ubf)(\mathbf{x}) \, d\mathbf{x}.
\end{equation}
For $|| \mathbf{h}||_{\mathcal{C}^1(\Omega)^d}$ small enough, we have $$
I (\ubf + \hbf) - I(\ubf) \geq (\boldsymbol{\xi}, \hbf),
$$
where $(\,,\,)$ is the duality product between $\mathcal{C}^1(\Omega)^d$ and its dual,
 which means that $\boldsymbol{\xi} \in \partial I$ with $\partial I$ being the subgradient of $I$.

\end{prop}
\begin{rem}
It is worth emphasizing that if $\varepsilon$ has no jump, then $I$ is Fr\'echet differentiable and $\boldsymbol{\xi}$ is its Fr\'echet derivative.
\end{rem}
\begin{rem}\label{defepstilde2} Under the assumptions of Proposition \ref{diffentiabilite},  if $\ubf^*$ is small enough (in $\mathcal{C}^1$-norm), then $\epstilde = \varepsilon \circ \left( \mathbb{I} + \ubf^*\right)^{-1} $ can be written as
\begin{equation}\label{formepstilde}
\epstilde=\epstilde_1  + \epstilde_2 \chi_{\tilde{\Omega}_2} ,
\end{equation}
with $\epstilde_1\in \mathcal{C}^1(\ol \Omega)$ and $\epstilde_2 \in \mathcal{C}^{1}_0(\Omega)$. In the sequel, we shall define $\tilde{\Omega}_i = \left(\mathbb{I} + \ubf^*\right)(\Omega_i)$ and $\tilde{f}_i = \varepsilon_i\circ \left(\mathbb{I} + \ubf^*\right)^{-1}$. For doing so, we extend $\tilde{f}_1$ into a function $\epstilde_1$ defined on the whole domain such that $\epstilde_1\in \mathcal{C}^1(\ol \Omega)$ and $\epstilde_1\big\vert_{\tilde{\Omega}_1}= \tilde{f}_1$. Then, we set $\epstilde_2 = \tilde{f}_2 - \epstilde_1$ on $\tilde{\Omega}_2$. Finally, we extend $\epstilde_2$ into a compactly supported $\mathcal{C}^1$-function on the whole domain $\Omega$.
\end{rem}
We first prove the following lemma.
\begin{lemm} \label{lastterm}Let $\ubf,\hbf \in \mathcal{C}^1(\Omega)^d$ and let $\epstilde$ be as in (\ref{formepstilde}).
Then, for $\Vert \ubf-\ubf^* \Vert_{\mathcal{C}^1(\Omega)^d}$ and $\Vert \hbf\Vert_{\mathcal{C}^1(\Omega)^d}$ small enough, we have
\begin{equation}
\int_{\Omega}\left[\epstilde(\mathbf{x} +\ubf + \hbf) -\epstilde(\mathbf{x}+\ubf )\right]^2\, d\mathbf{x}  = \int_\Omega  \epstilde_2^2(\mathbf{x}+\ubf) | \hbf \cdot \nubf |  \delta_{\dr \tilde{\Omega}_2}(\mathbf{x}+\ubf) \, d\mathbf{x}+ o(\Vert \hbf \Vert_{\mathcal{C}^1(\Omega)^d}),
\end{equation}
where $\delta_{\dr\tilde{\Omega}_2}$ is the Dirac distribution on
${\dr\tilde{\Omega}_2}$ and
$\epstilde_2$ is defined in Remark \ref{defepstilde2}.

\end{lemm}

\Proof We start by decomposing $\epstilde$ as follows:
\begin{multline*}
\int_{\Omega}\left[\epstilde(\mathbf{x}+\ubf + \hbf) -\epstilde(\mathbf{x}+\ubf )\right]^2\, d\mathbf{x}  = \\ \int_\Omega \bigg[ \big( \epstilde_1(\mathbf{x}+\ubf+\hbf) - \epstilde_1(\mathbf{x}+\ubf)\big) +\big(\epstilde_2 (\mathbf{x}+\ubf+\hbf) \chi_{\tilde{\Omega}_2 }(\mathbf{x}+\ubf+\hbf) - \epstilde_2(\mathbf{x}+\ubf) \chi_{\tilde{\Omega}_2}(\mathbf{x}+\ubf)\big)\bigg]^2\, d\mathbf{x}.
\end{multline*}
Now, by developing the square, the first term can be estimated by
\begin{equation*}
\left\vert \int_\Omega \big( \epstilde_1(\mathbf{x}+\ubf+\hbf) - \epstilde_1(\mathbf{x}+\ubf)\big)^2 \, d\mathbf{x} \right\vert \leq \Vert\epstilde_1\Vert^2_{\mathcal{C}^1(\Omega)} \Vert \hbf \Vert_{\mathcal{C}^1({\Omega})^d}^2.
\end{equation*}
Next, we write \begin{multline*}
\epstilde_2 (x+\ubf+\hbf) \chi_{\tilde{\Omega}_2 }(\mathbf{x}+\ubf+\hbf) - \epstilde_2(\mathbf{x}+\ubf) \chi_{\tilde{\Omega}_2}(\mathbf{x}+\ubf) =\left[ \epstilde_2(\mathbf{x}+\ubf+\hbf)-\epstilde_2(\mathbf{x}+\ubf)\right] \chi_{\tilde{\Omega}_2}(\mathbf{x}+\ubf + \hbf)\\ + \left[\chi_{\tilde{\Omega}_2}(\mathbf{x}+\ubf+\hbf) - \chi_{\tilde{\Omega}_2}(\mathbf{x}+\ubf)\right] \epstilde_2(\mathbf{x}+\ubf).
 \end{multline*} Since $\left(\epstilde_1(\mathbf{x}+\ubf+\hbf)-\epstilde_1(\mathbf{x}+\ubf) \right) \epstilde_2(\mathbf{x}+\ubf) \in \mathcal{C}^1_0(\Omega)$,  Proposition \ref{firstorder} yields
$$\begin{array}{l}
\displaystyle  \left\vert\int_\Omega\big[ \epstilde_1(\mathbf{x}+\ubf+\hbf) - \epstilde_1(\mathbf{x}+\ubf) \big] \left[\chi_{\tilde{\Omega}_2}(\mathbf{x}+\ubf+\hbf) - \chi_{\tilde{\Omega}_2}(\mathbf{x}+\ubf)\right] \epstilde_2(\mathbf{x}+\ubf) \, d\mathbf{x}\right\vert \\ \nm \displaystyle \leq  C \left( \int_\Omega \left[\hbf\cdot\g \epstilde_1(\mathbf{x}+\ubf) \right]^2 \, d\mathbf{x}\right)^{1/2} \left( \left[\int_\Omega \left[\hbf\cdot \nubf  \epstilde_2(\mathbf{x}+\ubf) \right]^2 \delta_{\partial \tilde{\Omega}_2}(\mathbf{x}+\ubf) \, d\mathbf{x} \right]+ o(\Vert \hbf\Vert_{\mathcal{C}^1({\Omega})^d})  \right)^{1/2}\\ \nm \displaystyle \leq C \Vert \hbf\Vert_{\mathcal{C}^1({\Omega})^d}^2.
 \end{array}$$
 We now need to handle the last term
 \begin{equation*}\begin{array}{l} \ds
 \int_\Omega \big(\left[\chi_{\tilde{\Omega}_2}(\mathbf{x}+\ubf+\hbf) - \chi_{\tilde{\Omega}_2}(
 \mathbf{x}+\ubf)\right] \epstilde_2(\mathbf{x} +\ubf)\big)^2 \, d\mathbf{x}\\ \nm \qquad \ds =  \int_\Omega \left\vert\chi_{\tilde{\Omega}_2}(\mathbf{x}+\ubf+\hbf) - \chi_{\tilde{\Omega}_2}(\mathbf{x}+\ubf)\right\vert \epstilde_2(\mathbf{x}+\ubf)^2\, d\mathbf{x}.
 \end{array}
 \end{equation*}
 Using Proposition \ref{firstorder}, we obtain that
 \begin{equation*}
 \int_\Omega \big(\left\vert\chi_{\tilde{\Omega}_2}(\mathbf{x}+\ubf+\hbf) - \chi_{\tilde{\Omega}_2}(\mathbf{x}+\ubf)\right\vert \epstilde_2(\mathbf{x}+\ubf)\big)^2 \, d\mathbf{x} = \int_\Omega  \epstilde_2^2(\mathbf{x}+\ubf) | \hbf \cdot \nubf |  \delta_{\dr \tilde{\Omega}_2}(\mathbf{x}+\ubf) \, d\mathbf{x}+ o(\Vert \hbf \Vert_{\mathcal{C}^1({\Omega})^d}),
 \end{equation*}
 which completes the proof of the lemma. \cqfd

We are now ready to prove Proposition \ref{diffentiabilite}.

\Proof If $\ubf \in \mathcal{C}^1(\Omega)^2 $ and $\hbf\in \mathcal{C}^1(\Omega)^2$, then we have
\begin{equation*}
I(\ubf +\hbf)-I(\ubf) = \int_{\Omega} \left[\epstilde(\mathbf{x} +\ubf + \hbf) + \epstilde(\mathbf{x}+\ubf )- 2 \varepsilon(\mathbf{x})\right] \left[\epstilde(\mathbf{x}+\ubf + \hbf) -\epstilde(\mathbf{x}+\ubf) \right] \, d\mathbf{x},
\end{equation*}
and hence,
\begin{multline*}
I(\ubf +\hbf)-I(\ubf) = \int_{\Omega} \left[\epstilde(\mathbf{x}+\ubf + \hbf) -\epstilde(\mathbf{x}+\ubf )\right]^2 \, d\mathbf{x} \\+ 2  \int_{\Omega} \left[\epstilde(\mathbf{x}+\ubf )-  \varepsilon(\mathbf{x})\right]  \left[\epstilde(\mathbf{x}+\ubf + \hbf) -\epstilde(\mathbf{x}+\ubf) \right] \, d\mathbf{x}.
\end{multline*}
For any $\eta>0$, let $g^{(\eta)}$ be a smooth, compactly supported function such that $$\Vert g^{(\eta)} -  \left[\epstilde\circ (\mathbb{I}+\ubf )-  \varepsilon\right] \Vert_{L^2(\Omega)} < \eta \quad \mbox{ and } \quad \big\vert \vert g^{(\eta)}\vert_{\mathrm{TV}(\Omega)} -\vert \epstilde\circ(\mathbb{I}+\ubf )-  \varepsilon\vert_{\mathrm{TV}(\Omega)}\big\vert <\eta; $$ see \cite{ambrosio2000functions}.

 Now, we write
 \begin{multline*}
 \int_{\Omega} \left[\epstilde(\mathbf{x}+\ubf )-  \varepsilon(\mathbf{x})\right]  \left[\epstilde(\mathbf{x}+\ubf + \hbf) -\epstilde(\mathbf{x}+\ubf) \right] \, d\mathbf{x} =   \int_{\Omega}  g^{\eta}(\mathbf{x}) \left[\epstilde(\mathbf{x}+\ubf + \hbf) -\epstilde(\mathbf{x}+\ubf) \right] \, d\mathbf{x} \\+\int_{\Omega}  \left[\epstilde(\mathbf{x}+\ubf )-  \varepsilon(\mathbf{x}) - g^{\eta}(\mathbf{x}) \right] \left[\epstilde(\mathbf{x}+\ubf + \hbf) -\epstilde(\mathbf{x}+\ubf) \right]\, d\mathbf{x}.
 \end{multline*}

Let $\tau_{\hbf}$ be the translation operator. Then, $\tau_\hbf$ satisfies, for any $\hbf\in \mathcal{C}^1(\Omega)^d$,
\begin{equation}\label{liptau}
\Vert \tau_\hbf[f]-f \Vert_p \leq C(f) \Vert \hbf \Vert_{\mathcal{C}^1({\Omega})^d}, \quad \forall \; f\in \mathrm{SBV}^p(\Omega).
\end{equation}
Using Cauchy-Schwartz' inequality,  we get
 \begin{equation}\label{eq1}
\bigg\vert \int_{\Omega}  \left[\epstilde(\mathbf{x}+\ubf )-  \varepsilon(\mathbf{x}) - g^{\eta}(\mathbf{x}) \right] \left[\epstilde(\mathbf{x}+\ubf + \hbf) -\epstilde(\mathbf{x}+\ubf) \right] \, d\mathbf{x} \bigg\vert \leq C \eta \Vert \hbf \Vert_{\mathcal{C}^1({\Omega})^d},
 \end{equation}
where $C$ is a constant depending on $\epstilde, \ \ubf,$ and $ \Omega$.

We know that for a certain  function $\rho$ such that $\rho(s)\rightarrow 0$ when $s \rightarrow 0$:
\begin{equation}\label{eq3}
\bigg\vert \int_{\Omega}  g^{\eta}(\mathbf{x}) \left[\epstilde(\mathbf{x}+\ubf + \hbf) -\epstilde(\mathbf{x}+\ubf) \right]
\, d\mathbf{x} - \int_{\Omega}  g^{\eta}(\mathbf{x}) \hbf(\mathbf{x}) \cdot D\left(\epstilde\circ(\mathbb{I}+\ubf)\right) \, d\mathbf{x} \bigg\vert \leq \Vert \hbf\Vert_{\mathcal{C}^1({\Omega})^d} \rho(\Vert\hbf\Vert_{\mathcal{C}^1({\Omega})^d}).
\end{equation}
Now, we have the following estimate:
\begin{equation}\label{eq2}
\bigg\vert \int_{\Omega}  g^{\eta}(\mathbf{x}) \hbf(\mathbf{x}) \cdot D\left(\epstilde\circ(\mathbb{I}+\ubf)\right) \, d\mathbf{x} - \int_{\Omega} [\epstilde(\mathbf{x}+\ubf) -\varepsilon(\mathbf{x})] \hbf(\mathbf{x})\cdot D\left(\epstilde\circ(\mathbb{I}+\ubf)\right) \, d\mathbf{x} \bigg\vert \leq C'\eta \Vert \hbf \Vert_{\mathcal{C}^1({\Omega})^d}.
\end{equation}

Indeed, since $\epstilde\in \piecewise \subset \mathrm{SBV}(\Omega)$, $\epstilde\circ(\mathbb{I}+\ubf)\in  \mathrm{SBV}(\Omega)$ and we can write the following decomposition of $D\left(\epstilde\circ(\mathbb{I}+\ubf)\right)$ into a continuous part and a jump part on a rectifiable surface $S$:
\begin{equation*}
D\left(\epstilde\circ(\mathbb{I}+\ubf)\right)=\g\left( \epstilde\circ(\mathbb{I}+\ubf)\right)\mathcal{H}^d+[\epstilde\circ(\mathbb{I}+\ubf)]\nubf_S \cH^{d-1}_S,
\end{equation*} we have that
\begin{equation*}
\bigg\vert \int_{\Omega}  \big[ g^{\eta}(\mathbf{x})-[\epstilde(\mathbf{x}+\ubf) -\varepsilon(\mathbf{x})]\big] \hbf(\mathbf{x}) \cdot\g\left( \epstilde\circ(\mathbb{I}+\ubf)\right)(\mathbf{x})\, d\mathbf{x}\bigg\vert \leq C_1\eta \Vert \hbf\Vert_{\mathcal{C}^1({\Omega})^d}.
\end{equation*}
For the jump part, since $S$ is a rectifiable surface and the function $f^\eta=g^{\eta}-[\epstilde\circ (\mathbb{I}+\ubf) -\varepsilon]$ is piecewise continuous, it is possible to define a trace   $f^\eta|_S$ on the surface $S$  satisfying $$\Vert f^\eta|_S \Vert_{L^1(S)} \leq C_2 \Vert f^\eta \Vert_{L^1(\Omega)}$$ for some positive constant $C_2$ depending only on $S$ and $\Omega$.
Then we get
\begin{equation*}
\bigg\vert \int_{S}  f^\eta \hbf(\mathbf{x}) \cdot[\epstilde\circ (\mathbb{I}+\ubf)]\nubf_S \cH^{d-1}_S \bigg\vert \leq C_3 \eta \Vert \hbf\Vert_{\mathcal{C}^1({\Omega})^d}
\end{equation*}
for some positive constant $C_3$ independent of $\eta$ and $\mathbf{h}$.

Now, the last term $ \int_{\Omega} \left[\epstilde(\mathbf{x}+\ubf + \hbf) -\epstilde(\mathbf{x} +\ubf )\right]^2$ can be handled using Lemma \ref{lastterm}. Doing so, we obtain
\begin{equation}\label{3terme}
 \int_{\Omega} \left[\epstilde(\mathbf{x} +\ubf + \hbf) -\epstilde(\mathbf{x} +\ubf )\right]^2 = \int_\Omega  \epstilde_2^2(\mathbf{x} +\ubf) | \hbf \cdot \nubf |  \delta_{\dr \tilde{\Omega}_2}(\mathbf{x}+\ubf) + o(\Vert \hbf \Vert_{\mathcal{C}^1({\Omega})^d}).
\end{equation}
Combining (\ref{eq1}), (\ref{eq3}), (\ref{eq2}), and  (\ref{3terme}), we get that for every $\eta >0$,
\begin{multline*}
\bigg\vert I(\ubf +\hbf)-I(\ubf) - 2\int_{\Omega} [\epstilde(\mathbf{x}+\ubf) -\varepsilon(\mathbf{x})] \hbf(\mathbf{x})\cdot D\epstilde\circ(\mathbb{I}+\ubf)(\mathbf{x}) \, d\mathbf{x} - \int_\Omega  \epstilde_2^2(\mathbf{x}+\ubf) | \hbf \cdot \nubf | \delta_{\dr \tilde{\Omega}_2}(\mathbf{x}+\ubf) \, d\mathbf{x} \bigg\vert \\ \leq
C_4  \Vert \hbf\Vert_{\mathcal{C}^1({\Omega})^d} \bigg( \rho(\Vert\hbf\Vert_{\mathcal{C}^1({\Omega})^d}) + \eta \bigg)
\end{multline*}
for some positive constant $C_4$ independent of $\mathbf{h}$ and $\eta$.

Finally, it follows that
$$
I (\ubf + \hbf) - I(\ubf) = (\boldsymbol{\xi}, \hbf) + \int_\Omega  \epstilde_2^2(\mathbf{x}+\ubf) | \hbf \cdot \nubf | \delta_{\dr \tilde{\Omega}_2}(\mathbf{x}+\ubf) \, d\mathbf{x} +
o(\Vert \hbf\Vert_{\mathcal{C}^1({\Omega})^d}),
$$
where $\boldsymbol{\xi}$ is defined by (\ref{defxi}). Hence, either $\ds \int_\Omega  \epstilde_2^2(\mathbf{x}+\ubf) | \hbf \cdot \nubf | \delta_{\dr \tilde{\Omega}_2}(\mathbf{x}+\ubf) \, d\mathbf{x}$ is of order of $\Vert \hbf\Vert_{\mathcal{C}^1({\Omega})^d}$ and we get
$$
I (\ubf + \hbf) - I(\ubf) \geq (\boldsymbol{\xi}, \hbf)
$$
for $\Vert \hbf\Vert_{\mathcal{C}^1({\Omega})^d}$ small enough or $\ds \int_\Omega  \epstilde_2^2(\mathbf{x}+\ubf) | \hbf \cdot \nubf | \delta_{\dr \tilde{\Omega}_2}(\mathbf{x}+\ubf) \, d\mathbf{x} = o(\Vert \hbf\Vert_{\mathcal{C}^1({\Omega})^d})$ and in this case, $I$ is Fr\'echet differentiable and $\boldsymbol{\xi}$ is its Fr\'echet derivative.
The proof of Proposition \ref{diffentiabilite} is then complete. \cqfd

\begin{rem} The minimization of the functional $I$ gives a reconstruction of $\ubf^*$ on a subdomain $\Omega \subset \Omega_0$. In practical conditions, since $\ubf^*$ is small $\Omega$ is almost the whole domain $\Omega_0$. The values of $\ubf^*$ on the boundary are known and, since $\ubf^*$ is of class $\mathcal{C}^1$, it is possible to deduce the values of $\ubf^*$ on $\Omega_0 \setminus \Omega$ by interpolation.
\end{rem}

\section{Reconstruction of the shear modulus} \label{sect4}
The problem is now to recover the function $\mu$ the reconstructed internal data $\ubf$. For doing so, we use the method described in \cite{ammarielasto}.
We introduce the operator $\mathcal{F}$
 \begin{equation*}\ubf=\mathcal{F}[\mu] =
\left\{\begin{aligned}
\g\cdot \left(\mu (\g \ubf + \g \ubf^T) \right)+ \g p=0 \quad &  \mbox{in } \ \Omega_0,\\
\g\cdot \ubf = 0 \quad & \mbox{in }\ \Omega_0,\\
\ubf = \mathbf{f} \quad & \mbox{on } \ \dr \Omega_0,
\end{aligned}\right.
\end{equation*}
and minimize the function $\mathcal{K}$ given by
 \begin{equation*}
\begin{aligned}
\mathcal{C}^{0,1}(\overline{\Omega_0}) &\longrightarrow  \mathbb{R} \\
\mu &\longmapsto\mathcal{K}[\mu]= \int_\Omega \vert \mathcal{F}[\mu]- \ubf \vert^2\, d\mathbf{x}.
\end{aligned}
\end{equation*}
According to \cite{ammarielasto}, $\mathcal{K}$ is Fr\'echet differentiable and its gradient can be explicitly computed. Let $\mathbf{v}$ be the solution of
\begin{equation*}
\left\{\begin{aligned}
\g\cdot \left(\mu (\g \mathbf{v} + \g \mathbf{v}^T) \right)+ \g q= \left( \mathcal{K}[\mu]-\ubf\right) \quad &  \mbox{in } \ \Omega_0,\\
\g\cdot \mathbf{v}= 0 \quad & \mbox{in } \ \Omega_0,\\
\mathbf{v} = 0 \quad & \mbox{on } \ \dr \Omega_0.
\end{aligned} \right.
\end{equation*}
Then, \begin{equation*}
\g \mathcal{K}(\mu) [h] = \int_{\Omega_0} h ( \g\mathbf{v} + \g \mathbf{v}^T): (\g \ubf + \g \ubf^T)\, d\mathbf{x}.
\end{equation*}
A gradient descent method can be applied in order to reconstruct $\mu$ from $\ubf$. We refer to \cite{ammarielasto} for more details.

\section{Numerical experiments} \label{sect5}

\begin{figure}[h]
\begin{center}
%
%
%
%
\begin{tikzpicture}[scale=0.6]

\begin{axis}[%
width=4.1in,
height=3.80333333333333in,
axis on top,
scale only axis,
xmin=-0.00167224080267559,
xmax=1.00167224080268,
ymin=-0.00167224080267559,
ymax=1.00167224080268,
colormap={mymap}{[1pt] rgb(0pt)=(0,0,1); rgb(63pt)=(0,1,0.5)},
colorbar,
point meta min=9.99223527007837,
point meta max=18.8134990530267
]
\addplot graphics [xmin=-0.0012531328320802,xmax=1.00125313283208,ymin=-0.0012531328320802,ymax=1.00125313283208] {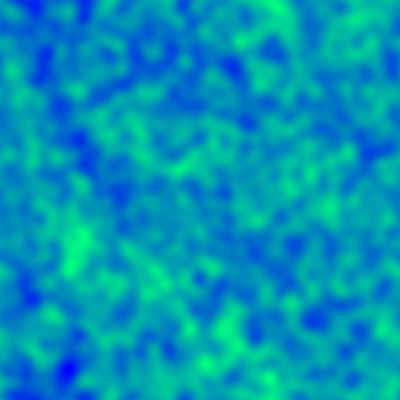};
\end{axis}
\end{tikzpicture}%
\caption{\label{epsilon} Optical index $\varepsilon$ of the medium.}
\end{center}
\end{figure}

\begin{figure}
\begin{center}
%
%
%
%
\begin{tikzpicture}[scale=0.6]

\begin{axis}[%
width=4.1in,
height=3.80333333333333in,
axis on top,
scale only axis,
xmin=-0.00167224080267559,
xmax=1.00167224080268,
ymin=-0.00167224080267559,
ymax=1.00167224080268,
colormap/jet,
colorbar,
title={Shear modulus distribution $\mu$},
point meta min=1.00000000824461,
point meta max=4.99982103508834
]

\addplot graphics [xmin=-0.00167224080267559,xmax=1.00167224080268,ymin=-0.00167224080267559,ymax=1.00167224080268] {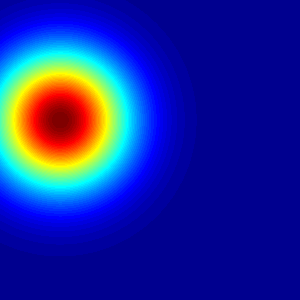};
\end{axis}
\end{tikzpicture}%
\begin{tikzpicture} [scale=0.6]

\begin{axis}[%
width=4.1in,
height=3.80333333333333in,
axis on top,
scale only axis,
xmin=-0.00167224080267559,
xmax=1.00167224080268,
ymin=-0.00167224080267559,
ymax=1.00167224080268,
colormap/jet,
colorbar,
title={Reconstructed shear modulus distribution $\mu_{\mathrm{rec}}$},
point meta min=0.5,
point meta max=4.1518157320979
]
\addplot graphics [xmin=-0.00167224080267559,xmax=1.00167224080268,ymin=-0.00167224080267559,ymax=1.00167224080268] {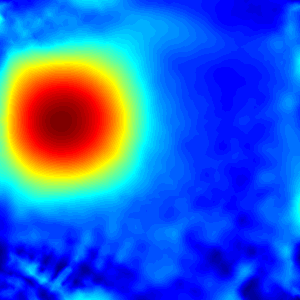};
\end{axis}
\end{tikzpicture}
\label{mu}
\caption{Shear modulus reconstruction.}
\end{center}
\end{figure}

We take $\Omega=[0,1]^2$ and discretize it on a $300\times 300$ grid, and generate a random Gaussian process to model the optical index $\varepsilon$ of the medium as shown in Figure \ref{epsilon}. Given a shear modulus $\mu$ map on $\Omega$; see Figure \ref{mu} (left), we solve  (\ref{eqmu}) on $\Omega$ via a finite element method compute the displacement field $\ubf$.
We then compute the displaced optical index $\varepsilon_u$ by using a spline interpolation approach  and proceed to recover the shear modulus
from the data $\varepsilon$ and $\varepsilon_u$ on the grid by the method described in the paper.

\begin{figure}
\begin{center}
%
%
%
%
\begin{tikzpicture}[scale=0.6]

\begin{axis}[%
width=4.1in,
height=3.80333333333333in,
axis on top,
scale only axis,
xmin=-0.00167224080267559,
xmax=1.00167224080268,
ymin=-0.00167224080267559,
ymax=1.00167224080268,
colormap/jet,
colorbar,
point meta min=0,
point meta max=0.99999999999981
]
\addplot graphics [xmin=-0.00167224080267559,xmax=1.00167224080268,ymin=-0.00167224080267559,ymax=1.00167224080268] {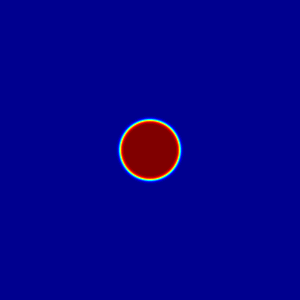};
\end{axis}
\end{tikzpicture}%
\caption{\label{ker}Averaging kernel $w_\delta$.}
\end{center}
\end{figure}

\begin{figure}
\begin{center}
%
%
%
%
\begin{tikzpicture}[scale=0.6]

\begin{axis}[%
width=4.1in,
height=3.80333333333333in,
axis on top,
scale only axis,
xmin=-0.00168350168350168,
xmax=1.0016835016835,
ymin=-0.00168350168350168,
ymax=1.0016835016835,
colormap/jet,
colorbar,
point meta min=1.00115231966457,
point meta max=2.90013207187996
]
\addplot graphics [xmin=-0.00168350168350168,xmax=1.0016835016835,ymin=-0.00168350168350168,ymax=1.0016835016835] {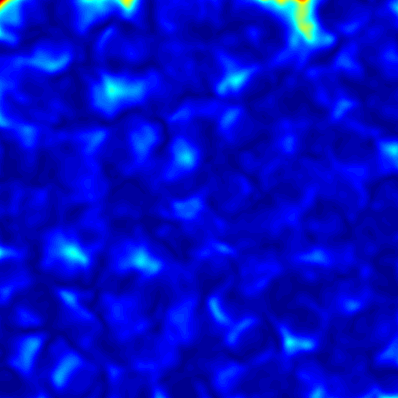};
\end{axis}
\end{tikzpicture}%
\caption{\label{condi} Conditioning of the matrix $w_\delta\star \g \varepsilon \g \varepsilon^T$.}
\end{center}
\end{figure}

\begin{figure}
\begin{center}
%
%
%
%

\begin{tikzpicture}

\begin{axis}[%
width=4.25cm,
height=3.25cm,
axis on top,
scale only axis,
xmin=0,
xmax=1,
ymin=0,
ymax=1,
name=plot1,
title={ $\ubf^* \cdot\mathbf{ e}_1$},
colormap/jet,
colorbar,
point meta min=0,
point meta max=0.005
]
\addplot graphics [xmin=-0,xmax=1,ymin=0,ymax=1]{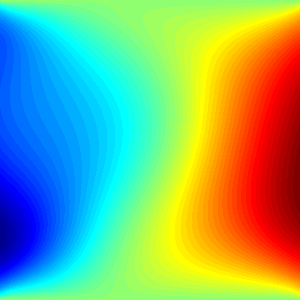};
\end{axis}

\end{tikzpicture}
\begin{tikzpicture}
\begin{axis}[%
width=4.25cm,
height=3.25cm,
axis on top,
scale only axis,
xmin=0,
xmax=1,
ymin=0,
ymax=1,
name=plot2,
title={ $\ubf^* \cdot\mathbf{ e}_2$},
colormap/jet,
colorbar,
point meta min=0,
point meta max=0.005
]
\addplot graphics [xmin=-0,xmax=1,ymin=0,ymax=1]{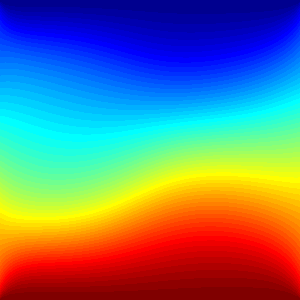};
\end{axis}

\end{tikzpicture}

\begin{tikzpicture}

\begin{axis}[%
width=4.25cm,
height=3.25cm,
axis on top,
scale only axis,
xmin=0,
xmax=1,
ymin=0,
ymax=1,
name=plot3,
title={Initial guess $\ubf_\delta\cdot \mathbf{e}_\text{1}$},
colormap/jet,
colorbar,
point meta min=0,
point meta max=0.005
]
\addplot graphics [xmin=-0,xmax=1,ymin=0,ymax=1] {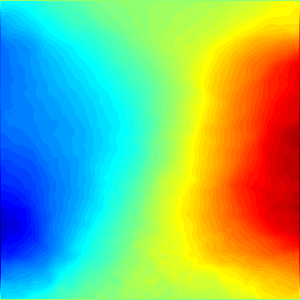};
\end{axis}
\end{tikzpicture}
\begin{tikzpicture}

\begin{axis}[%
width=4.25cm,
height=3.25cm,
axis on top,
scale only axis,
xmin=0,
xmax=1,
ymin=0,
ymax=1,
name=plot4,
title={Initial guess $\ubf_\delta\cdot \mathbf{e}_\text{2}$},
colormap/jet,
colorbar,
point meta min=0,
point meta max=0.005
]
\addplot graphics [xmin=-0,xmax=1,ymin=0,ymax=1]{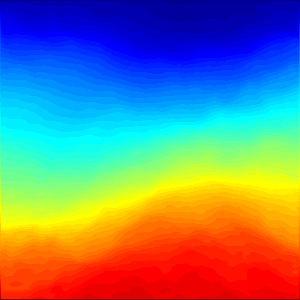};
\end{axis}
\end{tikzpicture}

\begin{tikzpicture}
\begin{axis}[%
width=4.25cm,
height=3.25cm,
axis on top,
scale only axis,
xmin=0,
xmax=1,
ymin=0,
ymax=1,
title={Reconstructed $\ubf\cdot\mathbf{e}_1$},
colormap/jet,
colorbar,
point meta min=0,
point meta max=0.005
]
\addplot graphics  [xmin=-0,xmax=1,ymin=0,ymax=1]{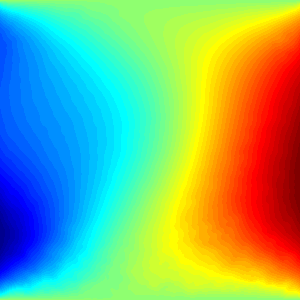};
\end{axis}
\end{tikzpicture}
\begin{tikzpicture}
\begin{axis}[%
width=4.25cm,
height=3.25cm,
axis on top,
scale only axis,
xmin=0,
xmax=1,
ymin=0,
ymax=1,
name=plot6,
title={Reconstructed $\ubf\cdot\mathbf{e}_2$},
colormap/jet,
colorbar,
point meta min=0,
point meta max=0.005
]
\addplot graphics  [xmin=-0,xmax=1,ymin=0,ymax=1]{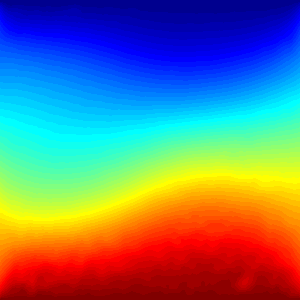};
\end{axis}
\end{tikzpicture}
\caption{\label{urec} Displacement field and its reconstruction.}
\end{center}
\end{figure}

Using (\ref{estinitial}), we first compute the initial guess $\ubf_\delta$ for the displacement field as the least-square solution to minimization of $J_{\mathbf{x}}$.
Figure \ref{ker} shows the kernel $w_\delta$ used to compute $\ubf_\delta$. As one can see $\delta$ needs to be large enough so the matrix $w_\delta \star \left( \g \varepsilon \g \varepsilon^T \right)$ is invertible at each point $\mathbf{x}$, which is basically saying that $\delta$ must be bigger than the correlation length of $\varepsilon$. Figure \ref{condi}  shows the conditioning of the matrix $w_\delta \star \left( \g \varepsilon \g \varepsilon^T \right)$.
Figure \ref{urec} shows the true displacement $\ubf^*$, the result of the first order  approximation (i.e., the initial guess) $\ubf_\delta$ and then the result of the optimization process using a gradient descent method to minimize the discrepancy functional $I$.

Once the displacement inside the domain is reconstructed, we can recover the shear modulus $\mu$, as shown in Figure \ref{mu}. We reconstruct $\mu$  by minimizing the
functional $\mathcal{K}$ and  using a gradient descent-type method. Note that
gradient of  $\mathcal{K}$ is computed with the adjoint state method,
described previously.  As it can be seen in Figure \ref{mu},
the reconstruction is very accurate  but not so perfect on the
boundaries of $\Omega$, which is due to the poor
estimation of $\ubf$ on $\partial \Omega$.

\section{Concluding remarks}

In this paper, we developed a novel algorithm which gives access not only to stiffness quantitative information of biological tissues but also opens the way to other contrasts such as mechanical anisotropy. In the heart, the muscle fibers have anisotropic mechanical properties. It would be very interesting to detect a change in fiber orientation using OCT elastographic tomography.

\end{document}